# ORDER REDUCTION OF NONLINEAR QUASI-PERIODIC SYSTEMS SUBJECTED TO EXTERNAL EXCITATIONS

**Sandesh G. Bhat**
Arizona State University
6075 S. Innovation Way West, 101D, Mesa, AZ 85212
sbhat9@asu.edu

**Susheelkumar Cherangara Subramanian**
Arizona State University
6075 S. Innovation Way West, 101, Mesa, AZ 85212
scherang@asu.edu

**Sangram Redkar**
Arizona State University
6075 S. Innovation Way West, 101B, Mesa, AZ 85212
sangram.redkar@asu.edu

## ABSTRACT

In his paper, we present order reduction techniques for nonlinear quasi-periodic systems subjected to external excitations. The order reduction techniques presented here are based on the Lyapunov-Perrone (L-P) Transformation. For a class of non-resonant quasi-periodic systems, the L-P transformation can convert a linear quasi-periodic system into a linear time-invariant one. This Linear Time-Invariant (LTI) system retains the dynamics of the original quasi-periodic system. Once this LTI system is obtained, the tools and techniques available for analysis of LTI systems can be used, and the results could be obtained for the original quasi-periodic system via the L-P transformation. This approach is similar to using the Lyapunov-Floquet (L-F) transformation to convert a linear time-periodic system into an LTI system and perform analysis and control.

Order reduction is a systematic way of constructing dynamical system models with relatively smaller states that accurately retain large-scale systems' essential dynamics. In this work, reduced-order modeling techniques for nonlinear quasi-periodic systems subjected to external excitations are presented. The methods proposed here use the L-P transformation that makes the linear part of transformed equations time-invariant. In this work, two order reduction techniques are suggested. The first method is simply an application of the well-known Guyan like reduction method to nonlinear systems. The second technique is based on the concept of an invariant manifold for quasi-periodic systems.

The 'quasi-periodic invariant manifold' based technique yields' reducibility conditions.' These conditions help us to understand the various types of resonant interactions in the system. These resonances indicate energy interactions between the system states, nonlinearity, and external excitation. To retain the essential dynamical characteristics, one has to preserve all these 'resonant' states in the reduced-order model. Thus, if the 'reducibility conditions' are satisfied then only, a nonlinear order reduction based on the quasi-periodic invariant manifold approach is possible. It is found that the invariant manifold approach yields good results. These methodologies are general and can be used for parametric study, sensitivity analysis, and controller design.

# INTRODUCTION

Order reduction is constructing small order systems from the large-scale structure, which captures the dominant dynamics [1]. In the design and development process, engineers are often faced with analyzing complex dynamical systems governed by a large set of integrodifferential (ordinary/partial) equations. These systems are complicated to solve analytically, and one has to resort to numerical techniques [2]. While solving these dynamical systems numerically, one has to consider various issues like convergence, numerical truncation errors, and, most importantly, the limited computational resources and time. To simulate the dynamical system's response accurately within a reasonable amount of time, one can try to construct an equivalent small-scale system known as the 'reduced order model' that will approximate the large-scale system dynamics [3,4]. This approach simplifies a sizeable dynamical system by replacing it with an equivalent small-scale system known as the 'Model Order Reduction' (MOR).

Researchers have used order reduction of linear systems using multiple techniques. Some of the techniques include error minimization [5], pole clustering [6], transfer function based [7], and Pade approximation [8], to mention a few. For a comprehensive overview, we refer to reference [9]. Nonlinear order reduction is studied by researchers from a structural point of view in the second-order and state-space form [10-12]. For order reduction of time-periodic system, researchers utilized the L-F transformation and performed order reduction. For details on this approach, we refer to references [13-17].

Order Reduction procedure comprises the following steps.
1. A study of the large-scale system and identifying the dominant states pertaining to the dominant dynamics.
2. Elimination of the non-dominant states either by simply neglecting their contribution or replacing them with the appropriate functions of dominant states.
3. Formulation of an equivalent reduced-order system consisting only of the dominant states.

In this work, we concentrate on the 'automatic' order reduction of an important class of engineering systems known as the 'nonlinear quasi-periodic systems.' [18-20]. These systems have quasi-periodic coefficients [7,8] and in the state-space form expressed as

$$\dot{\mathbf{x}} = \mathbf{A}(t)\mathbf{x} + \mathbf{f}(\mathbf{x},t) + \mathbf{F}(t) \tag{1}$$

where $\mathbf{A}(t)$ is an $n \times n$ matrix with quasi-periodic coefficients and $\mathbf{f}(\mathbf{x},t)$ is a nonlinear function with monomials of $\mathbf{x}$, $\mathbf{F}(t)$ is the external deterministic excitation, and $\mathbf{x}$ is an $n$ vector of appropriate dimensions. The objective of order reduction is to construct a reduced-order system

$$\dot{\mathbf{x}}_r = \mathbf{A}_r(t)\mathbf{x}_r + \mathbf{f}_r(\mathbf{x}_r,t) + \mathbf{F}_r(t) \tag{2}$$

that captures the essential dynamics of the large-scale system. It is noted that $\mathbf{A}(t)$ the matrix is quasi-periodic and contains incommensurate frequencies. To the best of the author's knowledge, no current techniques would allow direct order reduction from the equation (1) to the equation (2).

The reducibility of a linear quasi-periodic system has been the subject of research in the scientific community. Many excellent references discuss the reducibility of quasi-periodic systems [20-26]. Recently, Waswa and Redkar presented a technique based on L-F transformation, state augmentation, and normal form to reduce linear quasi-periodic system into an LTI system [21]. Very recently, Subramanian and Redkar presented a method to compute L-P transformation based on intuitive state augmentation and normal forms [27]. This technique can be utilized to calculate a closed-form expression for the L-P transformation. In this paper, we use the L-P transformation-based approach presented by Subramanian and Redkar to obtain the LTI representation of the quasi-periodic system and perform the order reduction. For clarity, it is noted that the reducibility of a quasi-periodic system means converting a linear quasi-periodic system into an LTI system of the same dimension, and order reduction means reducing the size (or the number of states) of the original quasi-periodic system or its LTI representation obtained via the L-P transformation.

This paper briefly reviews the L-P transformation computation and its application. The linear and nonlinear order reduction techniques are outlined in section two. In section three, we present an application- Mathieu Hill-type equation for which the L-P transformation can be determined using the state augmentation and normal forms. For this system, the reduced-order models are constructed using linear and nonlinear techniques. In the end, in section five, discussions and conclusions are presented.

## 2. MATHEMATICAL PRELIMINARIES

### 2.1 Computation of L-P transformation

A quasi-periodic dynamical system, without any external excitation, can be expressed as

$$\dot{\mathbf{y}} = \mathbf{A}(t)\mathbf{y}; \qquad (3)$$

where $\mathbf{A}(t)$ is a $n \times n$ matrix containing a finite number ($k$) of incommensurable frequencies ($k \geq 2$)

$$\mathbf{A}(t) = \tilde{A}(\omega_1 t, \ldots, \omega_k t) \quad \forall\, k \geq 2 \qquad (4)$$

It is noted that $\mathbf{A}(t)$ is continuous and periodic in each argument, but the ratio of any two frequencies is irrational[28]. It can be observed that the dynamical system given by the equation (4) is linear, and the method of normal forms [29] is inapplicable. However, the parametric excitation terms can be considered fictitious states, and the linear nonautonomous system given by the equation can be expressed as a nonlinear autonomous system as follows.

Consider the most general form of the equation (3) with quasi-periodicity provided by

$$\dot{\mathbf{x}} = (\mathbf{B}_0 + \mathbf{B}(t))\mathbf{x} \qquad (5)$$

where $\mathbf{A}(t)$ expressed as a constant matrix $\mathbf{B}_0$ and a quasi-periodic matrix $\mathbf{B}(t)$ of appropriate dimensions. Typically, $\mathbf{B}(t)$ comprises of $\sum_{i=1}^{n}(a_i \cos(\omega_i t) + b_i \sin(\omega_i t))$ type terms where $\omega_i$ is the frequency of quasi-periodic excitation (c.f. equation (4)). Assuming $q_i = \cos(\omega_i t)$ and $p_i = \sin(\omega_i t)$ the equation (5) can be expressed as

$$\dot{\bar{\mathbf{x}}} = \bar{\mathbf{B}}_0 \bar{\mathbf{x}} + \mathbf{f}(\bar{\mathbf{x}}) \qquad (6)$$

where $\bar{\mathbf{x}} = [\mathbf{x}, \mathbf{p}, \mathbf{q}]^T$, $\mathbf{p} = [p_1, p_2, \ldots p_n]^T$, $\mathbf{q} = [q_1, q_2, \ldots q_n]^T$ and

Applying the modal transformation $\bar{\mathbf{x}} = \mathbf{M}\mathbf{z}$, if $\mathbf{B}_0$ has semi-simple eigenvalues, the equation (6) is transformed into

$$\dot{\mathbf{z}} = \mathbf{J}\mathbf{z} + \mathbf{M}^{-1}\mathbf{f}(\mathbf{z}) \qquad (7)$$

where $\mathbf{J}$ is the Jordan form of $\bar{\mathbf{B}}_0$ (assumed to have semi-simple eigenvalues). The diagonal elements of the $\mathbf{J}$ matrix contain the linear matrix's ($\bar{\mathbf{B}}_0$) eigenvalues and frequencies of parametric excitations $[\lambda_{n+1}...\lambda_{n+m}]$ that are incommensurate.

$$\mathbf{J} = \begin{bmatrix} \lambda_1 & & & & \\ & \lambda_2 & & & \\ & & \ddots & & \\ & & & \lambda_n & \\ & & & & \omega \end{bmatrix}, \quad \omega = \begin{bmatrix} \lambda_{n+1} & & \\ & \ddots & \\ & & \lambda_{n+m} \end{bmatrix} \quad (8)$$

The system shown in the equation (7) is amenable to an application of Normal Forms[15,16]. A near identity transformation [30,31] (of the form given by equation (9)) is applied to the equation (7).

$$\mathbf{z} = \mathbf{v} + \mathbf{h}_r(\mathbf{v}) \quad (9)$$

where $\mathbf{h}_r(\mathbf{v})$ is a formal power series in $\mathbf{v}$ of degree $r$ with $T$ periodic coefficients that leads to

$$\dot{\mathbf{v}} = \mathbf{J}\mathbf{v} - \left[\frac{\partial \mathbf{h}_r(\mathbf{v})}{\partial(\mathbf{v})}\mathbf{J}\mathbf{v} - \mathbf{J}\mathbf{h}_r(\mathbf{v})\right] + \mathbf{f}_r(\mathbf{v}) \quad (10)$$

The higher-order nonlinear terms in the equation (10) are eliminated by considering the following condition

$$\frac{\partial \mathbf{h}_r(\mathbf{v})}{\partial(\mathbf{v})}\mathbf{J}\mathbf{v} - \mathbf{J}\mathbf{h}_r(\mathbf{v}) + \mathbf{f}_r(\mathbf{v}) = 0 \quad (11)$$

where

$$\mathbf{h}_r(\mathbf{v}) = \sum_{j=1}^{s}\sum_{\bar{m}_r}\sum_{v=-\alpha}^{v=\alpha} h_{j\bar{m}v}|\mathbf{v}|^m \mathbf{e}_j, \quad \mathbf{f}_r(\mathbf{v}) = \sum_{j=1}^{s}\sum_{\bar{m}_r}\sum_{v=-\alpha}^{v=\alpha} f_{j\bar{m}v}|\mathbf{v}|^m \mathbf{e}_j, \quad m_r = (m_1, m_2 ..), \sum_{i=1}^{n} m_i = 2, |\mathbf{v}|^m = v_1^{m_1} v_2^{m_2} .. v_n^{m_n}$$

and $\mathbf{e}_j$ is the $j^{th}$ member of the natural basis

After solving the equation (11), the solvability expression for a given degree of nonlinearity can be expressed as

$$h_{jmv} = \frac{f_{jmv}}{\mathbf{m}_r.\lambda - \lambda_j} \quad (12)$$

where $\lambda = [\lambda_1, \lambda_2, \lambda_3 .... \lambda_n, \lambda_{n+1} ... \lambda_{n+m}]^T$ are the eigenvalues of $\mathbf{J}$,

$$\mathbf{m}_r.\lambda - \lambda_j \neq 0 \quad (13)$$

In case the solvability condition [32] given by the equation (13) is satisfied, one can obtain the linear equation given by equation (14) and the near identify transformation given by equation (9).

$$\dot{\mathbf{v}} = \bar{\mathbf{J}}\mathbf{v} \quad (14)$$

The near identity transformation given by the equation (9) is wholly known in the non-resonant case and $\mathbf{h}_r(\mathbf{v})$ contains the terms that explicitly depend upon fictitious states $\mathbf{p}$ and $\mathbf{q}$. One can substitute these fictitious state in terms of their closed-form expression $q_i = \mathrm{Cos}(\omega_i t)$, and $p_i = \mathrm{Sin}(\omega_i t)$ that yields $\mathbf{h}_r(\mathbf{v},t)$ leading to the following form of the near identity transformation

$$\mathbf{z} = [\mathbf{I} + \overline{\mathbf{Q}}(t)]\mathbf{v} \approx \tilde{\mathbf{Q}}(t)\mathbf{v} \tag{15}$$

This transformation is similar to the Lyapunov-Floquet Transformation [33] but for quasi-periodic systems.

**Limitations:**

One important aspect is that this technique uses the normal form technique, which can be viewed as an extension of the higher-order averaging method [30,32] and has the same limitations as the averaging. This approach may not yield accurate results when the nonlinearity is very strong (i.e., very strong parametric excitation in the present case) or the linear term is absent (i.e., $\mathbf{B}_0 = \mathbf{0}$ c.f. equation (5)). On the other hand, to the best of the author's knowledge, this approach is the only approach that yields the L-P transformation given by the equation (15) in a closed-form. The authors have successfully used this approach to analyze linear and nonlinear quasi-periodic systems.

### 2.4 Computation of the inverse of the Lyapunov-Perrone Transformation:

For parametrically excited quasi-periodic linear systems of the equation (3) form, the L-P transformation is sufficient for carrying out analysis. The inverse of the L-P transformation is needed for quasi-periodic nonlinear systems or quasi-periodic linear/nonlinear systems with deterministic or stochastic excitations. The L-P transformation is a matrix where the matrix elements contain a truncated quasi-periodic Fourier series. Inverting a quasi-periodic matrix is not a trivial problem. In this section, we present two possible approaches to obtain the inverse L-P transformation.

**Symbolic Computation:** In minimal cases, when the L-P transformation matrix (given by equation (15)) is small ($2 \times 2$) and contains only a few terms, Symbolic computation software like Mathematica or Maple may be able to find the inverse. However, the inverse computed with this direct approach should be checked for the following conditions.

$$\begin{aligned} \tilde{\mathbf{Q}}^{-1}(0) &= \mathbf{I} \\ \tilde{\mathbf{Q}}^{-1}(t) \times \tilde{\mathbf{Q}}(t) &= \mathbf{I} \end{aligned} \tag{16}$$

The expression provided for $\mathbf{Q}^{-1}(t)$ may need further simplification for ease in future use.

**Neural Network:** One can also use a dynamical method using a recurrent neural network proposed for inversion of the time-varying matrix. One could use the gradient method [34], Zhang dynamics [35-37], or Chen dynamics [38] to find an inverse. In this section, we briefly present the Zhang dynamics approach[39] that could be used for inverting the L-P transformation.

Consider a time-varying matrix $\mathbf{Y}(t)$ with inverse $\mathbf{W}(t) = \mathbf{Y}^{-1}(t)$ so that the equation (17) is valid

$$\begin{aligned} \mathbf{Y}(t)\mathbf{W}(t) &= \mathbf{I} \\ \mathbf{Y}(t)\mathbf{W}(t) - \mathbf{I} &= \mathbf{0} \end{aligned} \tag{17}$$

We assume $\mathbf{Y}(t)$ is known and $\dfrac{d\mathbf{Y}(t)}{dt}$ exists. The objective is to find $\mathbf{W}(t)$ using the following equation

$$\mathbf{E}(\mathbf{W}(t), t) \equiv \mathbf{Y}(t)\mathbf{W}(t) - \mathbf{I} \tag{18}$$

where $\mathbf{E}(\mathbf{W}(t),t)$ is a matrix-valued error function. The derivative of the error function $\dot{\mathbf{E}}(\mathbf{W}(t),t)$ should be selected such that $\mathbf{E}(\mathbf{W}(t),t) \to \mathbf{0}$. Thus, $\dot{\mathbf{E}}(\mathbf{W}(t),t)$ can be chosen as

$$\frac{d\mathbf{E}(\mathbf{W}(t),t)}{dt} = -\Gamma \mathbf{F}(\mathbf{E}(\mathbf{W}(t),t)) \tag{19}$$

where $\Gamma$ is a scaling factor for the convergence and $\mathbf{F}(\mathbf{E}(\mathbf{W}(t),t))$ is called an activation function or matrix mapping recurrent neural network.

Differentiating equation (18) w.r.t. time and substituting equation (19) and (18) yields

$$\mathbf{Y}(t)\dot{\mathbf{W}}(t) = -\dot{\mathbf{Y}}(t)\mathbf{W}(t) - \Gamma \mathbf{F}(\mathbf{E}(\mathbf{W}(t),t))$$
$$\mathbf{Y}(t)\dot{\mathbf{W}}(t) = -\dot{\mathbf{Y}}(t)\mathbf{W}(t) - \Gamma \mathbf{F}(\mathbf{Y}(t)\mathbf{W}(t) - \mathbf{I}) \tag{20}$$

The equation (20) is a matrix differential equation that can be solved for $\mathbf{W}(t)$ using an appropriate initial condition. In the current paper $\mathbf{Y}(t)$ is the L-P transformation matrix $\tilde{\mathbf{Q}}(t)$ and $\mathbf{W}(t)$ is the inverse of L-P transformation $\tilde{\mathbf{Q}}^{-1}(t)$. Thus equation (20) can be written as

$$\tilde{\mathbf{Q}}(t)\dot{\tilde{\mathbf{Q}}}^{-1}(t) = -\dot{\tilde{\mathbf{Q}}}(t)\tilde{\mathbf{Q}}^{-1}(t) - \Gamma \mathbf{F}(\tilde{\mathbf{Q}}(t)\tilde{\mathbf{Q}}^{-1}(t) - \mathbf{I}) \tag{21}$$

One has to select an appropriate activation function and scaling constant $\Gamma \mathbf{F}()$ to achieve convergence. The equation (21) can be numerically integrated with the initial condition $\tilde{\mathbf{Q}}^{-1}(0) = \mathbf{I}$ to determine $\tilde{\mathbf{Q}}^{-1}(t)$. For more details on the Zhang Neural Network, its application, and proof of convergence, we refer to reference [37].

## 3. ORDER REDUCTION TECHNIQUES

### 3.1 Order reduction via linear projection

Consider a nonlinear quasi-periodic system described by the equation (1). Applying the L-P transformation $\mathbf{x}(t) = \mathbf{Q}(t)\mathbf{z}(t)$ produces

$$\dot{\mathbf{z}}(t) = \mathbf{J}\mathbf{z}(t) + \mathbf{Q}^{-1}(t)\mathbf{f}(\mathbf{z},t) + \mathbf{Q}^{-1}(t)\mathbf{F}(t) \equiv \mathbf{J}\mathbf{z}(t) + \mathbf{w}(\mathbf{z},t) + \bar{\mathbf{F}}(t) \tag{22}$$

where $\mathbf{J}$ is the constant matrix and $\mathbf{w}(\mathbf{z},t)$ represents an appropriately defined nonlinear quasi-periodic vector consisting of monomials of $z_j$.

Again, the objective of order reduction is to replace the nonlinear quasi-periodic system given by equation (22) with an equivalent system provided by

$$\dot{\mathbf{z}}_r(t) = \mathbf{J}_r \mathbf{z}_r(t) + \mathbf{w}_r(\mathbf{z}_r,t) + \mathbf{F}_r(t) \tag{23}$$

We partition equation (22) as

$$\begin{Bmatrix} \dot{\mathbf{z}}_r \\ \dot{\mathbf{z}}_s \end{Bmatrix} = \begin{bmatrix} \mathbf{J}_r & 0 \\ 0 & \mathbf{J}_s \end{bmatrix} \begin{Bmatrix} \mathbf{z}_r \\ \mathbf{z}_s \end{Bmatrix} + \begin{Bmatrix} \mathbf{w}_r(\mathbf{z}_r,\mathbf{z}_s,t) \\ \mathbf{w}_s(\mathbf{z}_r,\mathbf{z}_s,t) \end{Bmatrix} + \begin{Bmatrix} \mathbf{F}_r(t) \\ \mathbf{F}_s(t) \end{Bmatrix} \begin{matrix} (a) \\ (b) \end{matrix} \tag{24}$$

where $\mathbf{z}_s$ is an $(n-r)$ vector of non-dominant states, $\mathbf{J}_s$ is the matrix block of dimension $(n-r)\times(n-r)$ corresponding to the non-dominant states as defined earlier and $\mathbf{w}_r(\mathbf{z}_r,\mathbf{z}_s,t)$ and $\mathbf{w}_s(\mathbf{z}_r,\mathbf{z}_s,t)$ are the monomials of $\mathbf{z}$ (of order $i$) with quasi-periodic coefficients.

In the linear technique, the contribution of the non-dominant states is considered insignificant and hence neglected. Thus, the reduced-order model is given by

$$\dot{\mathbf{z}}_r(t) = \mathbf{J}_r \mathbf{z}_r(t) + \mathbf{w}_r(\mathbf{z}_r, 0, t) + \mathbf{F}_r(t) \tag{25}$$

The equation (25) is the reduced-order model of the actual large-scale system described by the equation (24). The equation (25) can be integrated numerically and using the transformation $\mathbf{x}(t) = \mathbf{Q}(t)\mathbf{T}\mathbf{z}_r(t)$, where $\mathbf{T} = \begin{bmatrix} \mathbf{I}_{r\times r} & \mathbf{0}_{r\times(n-r)} \end{bmatrix}^T$ all the states in $\mathbf{x}$ can be recovered.

This linear projection technique is simple and easy to implement. It may or may not provide accurate results. The selection of dominant states depends upon the judgment of the analyst. It does not give a clear insight into system dynamics if the system dynamics are complex and involves internal and parametric resonance.

**3.2 Order Reduction Using Invariant Manifold.**

This methodology is based on the ' Invariant Manifold Theory.' According to this theory, there exists a relationship between the dominant "master" and the non-dominant "slave" states of the system, and it is possible to replace (under certain conditions) the non-dominant states with dominant states. Thus, the order of the system can be reduced.

We assume that the frequency of forcing is incommensurate with the frequency of quasi-periodic parametric excitation. The constraint (or manifold governing) equations relating 'master' and 'slave' states are considerably complex but admit the solution in the form of asymptotic expansion. The relationship between the dominant and the non-dominant states of the system will involve contributions from the forcing and nonlinearity. If there are no resonances, then it is possible to replace the non-dominant states with the dominant states.

Once again, consider a nonlinear quasi-periodic system given by the equation (22) in the L-P transformed domain that is further partitioned as the equation (24). After ordering and expanding the nonlinear terms, we obtain.

$$\begin{aligned}
\dot{\mathbf{z}}_r &= \mathbf{J}_r \mathbf{z}_r + \varepsilon \mathbf{w}_{r2}(\mathbf{z}_r,\mathbf{z}_s,t) + \varepsilon^2 \mathbf{w}_{r3}(\mathbf{z}_r,\mathbf{z}_s,t) + \varepsilon^3 \mathbf{w}_{r4}(\mathbf{z}_r,\mathbf{z}_s,t) + ... \\
&\quad + \varepsilon^{i-1} \mathbf{w}_{ri}(\mathbf{z}_r,\mathbf{z}_s,t) + \varepsilon^{i-1} \mathbf{O}(|\mathbf{z}|^i) + \overline{\mathbf{F}}_r(t) \quad (a) \\
\dot{\mathbf{z}}_s &= \mathbf{J}_s \mathbf{z}_s + \varepsilon \mathbf{w}_{s2}(\mathbf{z}_r,\mathbf{z}_s,t) + \varepsilon^2 \mathbf{w}_{s3}(\mathbf{z}_r,\mathbf{z}_s,t) + \varepsilon^3 \mathbf{w}_{s4}(\mathbf{z}_r,\mathbf{z}_s,t) + ... \\
&\quad + \varepsilon^{i-1} \mathbf{w}_{si}(\mathbf{z}_r,\mathbf{z}_s,t) + \varepsilon^{i-1} \mathbf{O}(|\mathbf{z}|^i) + \overline{\mathbf{F}}_s(t) \quad (b)
\end{aligned} \tag{26}$$

where $\varepsilon^{n-1} \mathbf{w}_{rn}(\mathbf{z},t)$ include the terms of monomials of order $n$ in 'master' dynamics and $\varepsilon^{n-1} \mathbf{w}_{sn}(\mathbf{z},t)$ include terms of monomials of order $n$ in 'slave' dynamics. In this approach, we assume a nonlinear relationship between the dominant $(\mathbf{z}_r)$ and the non-dominant $(\mathbf{z}_s)$ states as

$$\mathbf{z}_s = \bar{\mathbf{h}}_{01}(t) + \varepsilon(\bar{\mathbf{h}}_{02}(t) + \bar{\mathbf{h}}_{12}(\mathbf{z}_r,t) + \bar{\mathbf{h}}_{22}(\mathbf{z}_r,t)) + \varepsilon^2(\bar{\mathbf{h}}_{03}(t) + \bar{\mathbf{h}}_{13}(\mathbf{z}_r,t) + \bar{\mathbf{h}}_{23}(\mathbf{z}_r,t) + \bar{\mathbf{h}}_{33}(\mathbf{z}_r,t))$$
$$+ \varepsilon^3(\bar{\mathbf{h}}_{04}(t) + \bar{\mathbf{h}}_{14}(\mathbf{z}_r,t) + \bar{\mathbf{h}}_{24}(\mathbf{z}_r,t) + \bar{\mathbf{h}}_{34}(\mathbf{z}_r,t) + \bar{\mathbf{h}}_{44}(\mathbf{z}_r,t)) + ...$$
(27)

Here $\bar{\mathbf{h}}_{ij}(\mathbf{z}_r,t)$ are the unknown quasi-periodic vector coefficients. Substitution of the equation (27) into (26) yields

$$\dot{\mathbf{z}}_s = \dot{\bar{\mathbf{h}}}_{01}(t) + \varepsilon\{\dot{\bar{\mathbf{h}}}_{02}(t) + \frac{\partial}{\partial t}(\bar{\mathbf{h}}_{12}(\mathbf{z}_r,t) + \bar{\mathbf{h}}_{22}(\mathbf{z}_r,t)) + \frac{\partial}{\partial \mathbf{z}_r}(\bar{\mathbf{h}}_{12}(\mathbf{z}_r,t) + \bar{\mathbf{h}}_{22}(\mathbf{z}_r,t)) \bullet \dot{\mathbf{z}}_r\}$$

$$+ \varepsilon^2\{\dot{\bar{\mathbf{h}}}_{03}(t) + \frac{\partial}{\partial t}(\bar{\mathbf{h}}_{13}(\mathbf{z}_r,t) + \bar{\mathbf{h}}_{23}(\mathbf{z}_r,t) + \bar{\mathbf{h}}_{33}(\mathbf{z}_r,t)) + \frac{\partial}{\partial \mathbf{z}_r}(\bar{\mathbf{h}}_{13}(\mathbf{z}_r,t) + \bar{\mathbf{h}}_{23}(\mathbf{z}_r,t) + \bar{\mathbf{h}}_{33}(\mathbf{z}_r,t)) \bullet \dot{\mathbf{z}}_r\}$$

$$+ \varepsilon^3\{\dot{\bar{\mathbf{h}}}_{04}(t) + \frac{\partial}{\partial t}(\bar{\mathbf{h}}_{14}(\mathbf{z}_r,t) + \bar{\mathbf{h}}_{24}(\mathbf{z}_r,t) + \bar{\mathbf{h}}_{34}(\mathbf{z}_r,t) + \bar{\mathbf{h}}_{44}(\mathbf{z}_r,t))$$

$$+ \frac{\partial}{\partial \mathbf{z}_r}(\bar{\mathbf{h}}_{14}(\mathbf{z}_r,t) + \bar{\mathbf{h}}_{24}(\mathbf{z}_r,t) + \bar{\mathbf{h}}_{34}(\mathbf{z}_r,t) + \bar{\mathbf{h}}_{44}(\mathbf{z}_r,t)) \bullet \dot{\mathbf{z}}_r\} + ...$$
(28)

Dropping spatial and temporal arguments for brevity, equation ((26)-b) can be rewritten as

$$\dot{\mathbf{z}}_s = \mathbf{J}_s \bullet \{\bar{\mathbf{h}}_{01} + \varepsilon(\bar{\mathbf{h}}_{02} + \bar{\mathbf{h}}_{12} + \bar{\mathbf{h}}_{22}) + \varepsilon^2(\bar{\mathbf{h}}_{03} + \bar{\mathbf{h}}_{13} + \bar{\mathbf{h}}_{23} + \bar{\mathbf{h}}_{33}) + \varepsilon^3(\bar{\mathbf{h}}_{04} + \bar{\mathbf{h}}_{14} + \bar{\mathbf{h}}_{24} + \bar{\mathbf{h}}_{34} + \bar{\mathbf{h}}_{44}) + ...$$

$$\sum_{m=4}^{n} \varepsilon^m (\mathbf{h}_{0m+1} + \sum_{k=1}^{m+1} \mathbf{h}_{k\,m+1}), t)$$

$$+ \varepsilon \mathbf{w}_{s2}(\bar{\mathbf{h}}_{01} + \varepsilon(\bar{\mathbf{h}}_{02} + \bar{\mathbf{h}}_{12} + \bar{\mathbf{h}}_{22}) + \varepsilon^2(\bar{\mathbf{h}}_{03} + \bar{\mathbf{h}}_{13} + \bar{\mathbf{h}}_{23} + \bar{\mathbf{h}}_{33}) + \varepsilon^3(\bar{\mathbf{h}}_{04} + \bar{\mathbf{h}}_{14} + \bar{\mathbf{h}}_{24} + \bar{\mathbf{h}}_{34} + \bar{\mathbf{h}}_{44}) + ...$$

$$\sum_{m=4}^{n} \varepsilon^m (\mathbf{h}_{0m+1} + \sum_{k=1}^{m+2} \mathbf{h}_{k\,m+1}), t)$$

$$+ \varepsilon^2 \mathbf{w}_{s3}(\bar{\mathbf{h}}_{01} + \varepsilon(\bar{\mathbf{h}}_{02} + \bar{\mathbf{h}}_{12} + \bar{\mathbf{h}}_{22}) + \varepsilon^2(\bar{\mathbf{h}}_{03} + \bar{\mathbf{h}}_{13} + \bar{\mathbf{h}}_{23} + \bar{\mathbf{h}}_{33}) + \varepsilon^3(\bar{\mathbf{h}}_{04} + \bar{\mathbf{h}}_{14} + \bar{\mathbf{h}}_{24} + \bar{\mathbf{h}}_{34} + \bar{\mathbf{h}}_{44}) + ...$$

$$\sum_{m=4}^{n} \varepsilon^m (\mathbf{h}_{0m+1} + \sum_{k=1}^{m+1} \mathbf{h}_{k\,m+1}), t)$$

$$+ \varepsilon^3 \mathbf{w}_{s4}(\bar{\mathbf{h}}_{01} + \varepsilon(\bar{\mathbf{h}}_{02} + \bar{\mathbf{h}}_{12} + \bar{\mathbf{h}}_{22}) + \varepsilon^2(\bar{\mathbf{h}}_{03} + \bar{\mathbf{h}}_{13} + \bar{\mathbf{h}}_{23} + \bar{\mathbf{h}}_{33}) + \varepsilon^3(\bar{\mathbf{h}}_{04} + \bar{\mathbf{h}}_{14} + \bar{\mathbf{h}}_{24} + \bar{\mathbf{h}}_{34} + \bar{\mathbf{h}}_{44}) + ...$$

$$\sum_{m=4}^{n} \varepsilon^m (\mathbf{h}_{0m+1} + \sum_{k=1}^{m+1} \mathbf{h}_{k\,m+1}), t)$$

$$+ \varepsilon^{i-1} \mathbf{w}_{si}(\bar{\mathbf{h}}_{01} + \varepsilon(\bar{\mathbf{h}}_{02} + \bar{\mathbf{h}}_{12} + \bar{\mathbf{h}}_{22}) + \varepsilon^2(\bar{\mathbf{h}}_{03} + \bar{\mathbf{h}}_{13} + \bar{\mathbf{h}}_{23} + \bar{\mathbf{h}}_{33}) + \varepsilon^3(\bar{\mathbf{h}}_{04} + \bar{\mathbf{h}}_{14} + \bar{\mathbf{h}}_{24} + \bar{\mathbf{h}}_{34} + \bar{\mathbf{h}}_{44}) + ...$$

$$\sum_{m=4}^{n} \varepsilon^m (\mathbf{h}_{0m+1} + \sum_{k=1}^{m+1} \mathbf{h}_{k\,m+1}), t))\}$$

$$+ \bar{\mathbf{F}}_s(t)$$
(29)

Substituting equation (26) and (27) into the equation (28) and equating it to the equation (29) yields a complex partial differential equation involving various orders of $\varepsilon$. By correlating the terms of the same order of $\varepsilon$, we obtain the equations, which need to be solved to determine $\mathbf{h}_{0m+1}(t)$ and $\mathbf{h}_{k\,m+1}(\mathbf{z}_m,t)$

Collecting the terms in order of $\varepsilon^0$ yields

$$\dot{\mathbf{h}}_{01}(t) = \mathbf{J}_s \bar{\mathbf{h}}_{01}(t) + \bar{\mathbf{F}}_s(t) \tag{30}$$

The equation (30) is a linear equation involving pure temporal arguments. The solution of the equation (30) can be determined using the convolution integral [29,40] as

$$\bar{\mathbf{h}}_{01}(t) = \mathbf{e}^{\mathbf{J}_s t}\mathbf{h}_{01}(0) + \int_0^t \mathbf{e}^{\mathbf{J}_s(t-\tau)} \bar{\mathbf{F}}_s(\tau) d\tau \tag{31}$$

If the forcing $\mathbf{F}(t)$ is harmonic with frequency $k\omega_f$, then after the L-P transformation, the frequency of harmonic excitation $\bar{\mathbf{F}}(t)$ becomes $\sum_{p_1=-\infty}^{+\infty} \sum_{p_2=-\infty}^{\infty} (\bar{\mathbf{p}} \bullet \boldsymbol{\omega}_p + k\omega_f)$, where $\boldsymbol{\omega}_p$ is the vector containing quasi-periodic frequencies in the L-P transformation $\boldsymbol{\omega}_p = \{\omega_1 \quad \omega_2\}, \bar{\mathbf{p}} = \{p_1 \quad p_2\}^T$.

Expressing forcing in the most general form as

$$\bar{\mathbf{F}}_s(t) = \sum_{k=-\infty}^{+\infty} \sum_{p_1=-\infty}^{+\infty} \sum_{p_2=-\infty}^{\infty} \mathbf{C}_{p_1 p_2 k} e^{\bar{i}(\mathbf{p} \bullet \boldsymbol{\omega}_p + k\omega_f)t} \tag{32}$$

If the eigenvalues of $\mathbf{J}_s$ are purely imaginary and given by $\bar{\lambda}_p$; $p = 1, 2, \ldots s$ then the solution can be written as

$$\bar{\mathbf{h}}_{01}(t) = \sum_{j=1}^{s} \sum_{p_1=-\infty}^{\infty} \sum_{p_2=-\infty}^{\infty} \sum_{k=-\infty}^{\infty} (C_{jp_1 p_2 k} \frac{e^{\bar{i}(\bar{n}\omega_p + k\omega_f)t}}{\bar{i}(\mathbf{p} \bullet \boldsymbol{\omega}_p + k\omega_f - \lambda_j)} e_j - C_{jp_1 p_2 k} \frac{e^{\bar{i}\lambda_j t}}{\bar{i}(\mathbf{p} \bullet \boldsymbol{\omega}_p + k\omega_f - \lambda_j)} e_j) \tag{33}$$

It can be seen that if $\mathbf{p} \bullet \boldsymbol{\omega}_p + k\omega_f - \lambda_j = 0$ for any combination, $\bar{\mathbf{h}}_{01}(t)$ cannot be found out, and the system is said to be in 'linear resonance.' This resonance is referred to as a 'primary' or a 'main resonance' in perturbation analysis.

Collecting the terms at the order of $\varepsilon^1$ yields

$$\dot{\mathbf{h}}_{02}(t) = \mathbf{J}_s \bar{\mathbf{h}}_{02}(t) - \frac{\partial \bar{\mathbf{h}}_{12}}{\partial \mathbf{z}_r} \mathbf{F}_r + \mathbf{w}_{s 2_0}(t) \tag{34}$$

$$\frac{\partial \bar{\mathbf{h}}_{12}}{\partial t} + \frac{\partial \bar{\mathbf{h}}_{12}}{\partial \mathbf{z}_r} \mathbf{J}_r \mathbf{z}_r + \frac{\partial \bar{\mathbf{h}}_{22}}{\partial \mathbf{z}_r} \mathbf{F}_r(t) - \mathbf{J}_s \bar{\mathbf{h}}_{12} = \mathbf{w}_{s 2_1}(\mathbf{z}_r, t) \tag{35}$$

$$\frac{\partial \bar{\mathbf{h}}_{22}}{\partial t} + \frac{\partial \bar{\mathbf{h}}_{22}}{\partial \mathbf{z}_r} \mathbf{J}_r \mathbf{z}_r - \mathbf{J}_s \bar{\mathbf{h}}_{22} = \mathbf{w}_{s 2_2}(\mathbf{z}_r, t) \tag{36}$$

It can be seen that equations (34), (35) and (36) are coupled equations. However, the equation (36) can be solved independently. Assuming the most general form of nonlinearity and expanding the known and unknown terms in multiple Fourier series as

$$\bar{\mathbf{h}}_{22}(\mathbf{z}_r,t) = \sum_{j=1}^{s} \sum_{\bar{\mathbf{m}}} \sum_{p_1=-\infty}^{+\infty} \sum_{p_2=-\infty}^{+\infty} h_{j\bar{\mathbf{m}} p_1 p_2 v}\, e^{\bar{i}(\bar{\mathbf{p}}\bullet\boldsymbol{\omega})t} \left|\mathbf{z}_r\right|^{\mathbf{m}} e_j \qquad (37)$$

$$\mathbf{w}_{s2_2}(\mathbf{z}_r,t) = \sum_{j=1}^{s} \sum_{\bar{\mathbf{m}}} \sum_{p_1=-\infty}^{+\infty} \sum_{p_2=-\infty}^{\infty} a_{j\bar{\mathbf{m}} p_1 p_2 v}\, e^{\bar{i}(\bar{\mathbf{p}}\bullet\boldsymbol{\omega})t} \left|\mathbf{z}_r\right|^{\mathbf{m}} e_j \qquad (38)$$

$$\left|\mathbf{z}_r\right|^{\mathbf{m}} = z_1^{m_1} z_2^{m_2} ... z_r^{m_r},\ \bar{i}=\sqrt{-1},\ m_1+......+m_r = 2;\ \boldsymbol{\omega}_p = \{\omega_1\ \ \omega_2\},\ \bar{\mathbf{p}} = \{p_1\ \ p_2\}^T$$

Collecting the terms to solve for the unknowns yields the '*reducibility condition*' given by the equation (39),

$$\bar{i}(\bar{\mathbf{p}}\bullet\boldsymbol{\omega}) + \sum_{l=1}^{r}(m_l \lambda_l) - \bar{\lambda}_p \neq 0 \qquad (39)$$

where all the terms appearing in the equation (39) are defined before.

In the absence of resonances $\bar{\mathbf{h}}_{22}(\mathbf{z}_r,t)$ can be obtained. It should be noted at this stage that forcing frequency does not appear in the equation (39), implying that there is no direct interaction between the nonlinearity and the external excitation. However, as we construct the solution using equations (34) and (35) forcing interacts with the nonlinearity, giving rise to additional 'resonance conditions.'

To find out the solution to the equation (35), which contains the contribution from the nonlinearity $\mathbf{w}_{s2}(\mathbf{z}_r,t)$, we expand the known and the unknown terms in the multiple Fourier series of the form

$$\bar{\mathbf{h}}_{12}(\mathbf{z}_r,t) = \sum_{j=1}^{s} \sum_{\bar{\mathbf{m}}} \sum_{p_1=-\infty}^{+\infty} \sum_{p_2=-\infty}^{+\infty} \sum_{p_3=-\infty}^{+\infty} h_{j\bar{\mathbf{m}} p_1 p_2 v}\, e^{\bar{i}(\tilde{\mathbf{p}}\bullet\bar{\boldsymbol{\omega}})t} \left|\mathbf{z}_r\right|^{\mathbf{m}} e_j \qquad (40)$$

$$\mathbf{w}_{rs_1}(\mathbf{z}_r,t) = \sum_{j=1}^{s} \sum_{\bar{\mathbf{m}}} \sum_{p_1=-\infty}^{+\infty} \sum_{p_2=-\infty}^{\infty} \sum_{p_3=-\infty}^{+\infty} a_{j\bar{\mathbf{m}} p_1 p_2 v}\, e^{\bar{i}(\tilde{\mathbf{p}}\bullet\bar{\boldsymbol{\omega}})t} \left|\mathbf{z}_r\right|^{\mathbf{m}} e_j \qquad (41)$$

$$\left|\mathbf{z}_r\right|^{\mathbf{m}} = z_1^{m_1} z_2^{m_2} ... z_r^{m_r},\ \bar{i}=\sqrt{-1},\ m_1+......+m_r = 1;\ \bar{\boldsymbol{\omega}} = \{\omega_1\ \ \omega_2\ \ \omega_f\};\ \tilde{\mathbf{p}} = \{p_1\ \ p_2\ \ p_3\}^T$$

A term-by-term comparison yields

$$h_{j\bar{m}v} = \frac{a_{j\bar{m}v}}{\bar{i}(\tilde{\mathbf{p}}\bullet\bar{\boldsymbol{\omega}}) + \lambda_l - \bar{\lambda}_p} \qquad (42)$$

The 'combined reducibility condition' can be expressed as

$$\bar{i}(\tilde{\mathbf{p}}\bullet\bar{\boldsymbol{\omega}}) + \lambda_l - \bar{\lambda}_p \neq 0 \qquad (43)$$

It can be observed that all the terms appearing in the equation (34) are free from spatial arguments, and it can be solved using the convolution, as discussed before. The forcing terms $\mathbf{w}_{s2_0}(t)$ (with a square type term in $\bar{\mathbf{h}}_{01}(t)$) appearing due to $\mathbf{w}_{s2}(\mathbf{z}_r, t)$ is known from the equation (31) which can be expressed in the form

$$\mathbf{w}_{s2_0}(t) = \sum_{j=1}^{s} \sum_{p_1=-\infty}^{+\infty} \sum_{p_2=-\infty}^{\infty} \sum_{p_3=-\infty}^{+\infty} a_{j\bar{m}v} e^{\bar{i}(\tilde{\mathbf{p}} \bullet \bar{\boldsymbol{\omega}})t} \, e_j \tag{44}$$

where $\bar{i} = \sqrt{-1}$, $\bar{\boldsymbol{\omega}} = \{\omega_1 \quad \omega_2 \quad \omega_f\}$; $\tilde{\mathbf{p}} = \{p_1 \quad p_2 \quad p_3\}^T$ and $\dfrac{\partial \bar{\mathbf{h}}_{12}}{\partial \mathbf{z}_r} \bar{\mathbf{F}}_r$ is known from the equation (33).

However, the equation (34) cannot be solved if $\lambda_s = \tilde{\mathbf{p}} \bullet \bar{\boldsymbol{\omega}}$, which can be written as

$$\lambda_s = p_1 \omega_p + p_2 \omega_{nl} + p_3 \omega_f \tag{45}$$

The exact combination will be determined by the kind of terms present in the forcing.

As we collect the terms at the order of $\varepsilon^2$, we obtain

$$\dot{\bar{\mathbf{h}}}_{03}(t) = \mathbf{J}_s \bar{\mathbf{h}}_{03}(t) + \frac{\partial \bar{\mathbf{h}}_{13}}{\partial \mathbf{z}_r} \bar{\mathbf{F}}_r(t) - \frac{\partial \bar{\mathbf{h}}_{12}}{\partial \mathbf{z}_r} \underbrace{\mathbf{w}_{m2_0}(t)}_{[\bar{\mathbf{h}}_{01}(t)]_2} + \underbrace{\mathbf{w}_{s2_{30}}(t)}_{\bar{\mathbf{h}}_{01}(t) \bullet \bar{\mathbf{h}}_{02}(t)} + \underbrace{\mathbf{w}_{s3_0}(t)}_{[\bar{\mathbf{h}}_{01}(t)]_3} \tag{46}$$

$$\frac{\partial \bar{\mathbf{h}}_{13}}{\partial t} + \frac{\partial \bar{\mathbf{h}}_{22}}{\partial \mathbf{z}_r} \underbrace{\mathbf{w}_{m2_0}(\mathbf{z}_r, t)}_{[\bar{\mathbf{h}}_{01}(t)]_2} + \frac{\partial \bar{\mathbf{h}}_{13}}{\partial \mathbf{z}_r} \mathbf{J}_r \mathbf{z}_r + \frac{\partial \bar{\mathbf{h}}_{23}}{\partial \mathbf{z}_r} \bar{\mathbf{F}}_r = \mathbf{J}_s \bar{\mathbf{h}}_{13} + \underbrace{\mathbf{w}_{s2_{31}}(\mathbf{z}_r, t)}_{\bar{\mathbf{h}}_{02}(t) \bullet \mathbf{z}_r + \bar{\mathbf{h}}_{01}(t) \bullet \bar{\mathbf{h}}_{12}(\mathbf{z}_r, t)} + \underbrace{\mathbf{w}_{s3_1}(\mathbf{z}_r, t)}_{[\bar{\mathbf{h}}_{01}(t)]_2 \bullet \mathbf{z}_r} \tag{47}$$

$$\frac{\partial \bar{\mathbf{h}}_{23}}{\partial t} + \frac{\partial \bar{\mathbf{h}}_{23}}{\partial \mathbf{z}_r} \mathbf{J}_r \mathbf{z}_r + \frac{\partial \bar{\mathbf{h}}_{33}}{\partial \mathbf{z}_r} \bar{\mathbf{F}}_r = \mathbf{J}_s \bar{\mathbf{h}}_{23} + \underbrace{\mathbf{w}_{s2_{32}}(\mathbf{z}_r, t)}_{\bar{\mathbf{h}}_{12}(\mathbf{z}_r, t) \bullet \mathbf{z}_r + \bar{\mathbf{h}}_{01}(t) \bullet \bar{\mathbf{h}}_{22}(\mathbf{z}_r, t)} + \underbrace{\mathbf{w}_{s3_2}(\mathbf{z}_r, t)}_{\bar{\mathbf{h}}_{01}(t) \bullet \mathbf{z}_r} \tag{48}$$

$$\frac{\partial \bar{\mathbf{h}}_{33}}{\partial t} + \frac{\partial \bar{\mathbf{h}}_{33}}{\partial \mathbf{z}_r} \mathbf{J}_r \mathbf{z}_r = \mathbf{J}_s \bar{\mathbf{h}}_{33} + \underbrace{\mathbf{w}_{s2_{33}}(\mathbf{z}_r, t)}_{\bar{\mathbf{h}}_{22}(\mathbf{z}_r, t) \bullet \mathbf{z}_r} + \mathbf{w}_{s3_3}(\mathbf{z}_r, t) \tag{49}$$

The equation (46) has only temporal arguments; the equation (47) is linear in spatial arguments; the equation (48) depends upon quadratic spatial arguments, and the equation (49) involves cubic spatial arguments. As before, one has to solve these equations sequentially.

It can be observed that the equation (49) can be solved independently and involves contribution from $\mathbf{w}_{s2}(\mathbf{z}_r, t)$ denoted by $\mathbf{w}_{s2_{33}}(\mathbf{z}_r, t)$. To solve this equation, we expand the known terms and unknown terms in the multiple Fourier series (c.f. equation(36))

$$\bar{\mathbf{h}}_{33}(\mathbf{z}_r,t) = \sum_{j=1}^{S}\sum_{\bar{\mathbf{m}}}\sum_{p_1=-\infty}^{+\infty}\sum_{p_2=-\infty}^{+\infty} h_{j\bar{m}v}\, e^{\bar{i}(\bar{\mathbf{p}}\bullet\boldsymbol{\omega})t}\,|\mathbf{z}_r|^{\mathbf{m}}\, e_j \tag{50}$$

$$\mathbf{w}_{s3_2}(\mathbf{z}_r,t) = \sum_{j=1}^{S}\sum_{\bar{\mathbf{m}}}\sum_{p_1=-\infty}^{+\infty}\sum_{p_2=-\infty}^{\infty} a_{j\bar{m}v}\, e^{\bar{i}(\bar{\mathbf{p}}\bullet\boldsymbol{\omega})t}\,|\mathbf{z}_r|^{\mathbf{m}}\, e_j \tag{51}$$

$$|\mathbf{z}_r|^{\mathbf{m}} = z_1^{m_1} z_2^{m_2}...z_r^{m_r},\ \bar{i}=\sqrt{-1},\ m_1+......+m_r = 3;\ \boldsymbol{\omega} = \{\omega_1\ \ \omega_2\},\ \bar{\mathbf{p}}=\{p_1\ \ p_2\}^T$$

It is possible to determine $\bar{\mathbf{h}}_{33}(\mathbf{z}_r,t)$ if the following 'reducibility condition' is satisfied.

$$\bar{i}\,(\bar{\mathbf{p}}\bullet\boldsymbol{\omega}) + \sum_{l=1}^{r}(m_l\lambda_l) - \bar{\lambda}_p \neq 0 \tag{52}$$

Once $\bar{\mathbf{h}}_{33}(\mathbf{z}_r,t)$ is known, we can solve the equation (48). This equation contains terms arising from the product of $\dfrac{\partial \bar{\mathbf{h}}_{33}(\mathbf{z}_r,t)}{\partial \mathbf{z}_r}\bar{\mathbf{F}}_r(t)$ (where $\bar{\mathbf{F}}_r(t)$ is the forcing on the master states), the contribution from quadratic $\mathbf{w}_{s2}(\mathbf{z}_r,t)$ nonlinearity (denoted by $\mathbf{w}_{s2_{32}}(\mathbf{z}_r,t)$), and contribution from cubic nonlinearity $\mathbf{w}_{s3}(\mathbf{z}_r,t)$ (represented by $\mathbf{w}_{s3_2}(\mathbf{z}_r,t)$). As before, we expand the known terms (marked as $\mathbf{w}_{sk}(\mathbf{z}_r,t)$) and the unknown terms ($\bar{\mathbf{h}}_{23}(\mathbf{z}_r,t)$) in multiple Fourier series of the form

$$\bar{\mathbf{h}}_{23}(\mathbf{z}_r,t) = \sum_{j=1}^{S}\sum_{\bar{\mathbf{m}}}\sum_{p_1=-\infty}^{+\infty}\sum_{p_2=-\infty}^{+\infty}\sum_{p_3=-\infty}^{+\infty} h_{j\bar{m}p_1p_2p_3v}\, e^{\bar{i}(\tilde{\mathbf{p}}\bullet\bar{\boldsymbol{\omega}})t}\,|\mathbf{z}_r|^{\mathbf{m}}\, e_j \tag{53}$$

$$\mathbf{w}_{sk}(\mathbf{z}_r,t) = \sum_{j=1}^{S}\sum_{\bar{\mathbf{m}}}\sum_{p_1=-\infty}^{+\infty}\sum_{p_2=-\infty}^{\infty}\sum_{p_3=-\infty}^{+\infty} a_{j\bar{m}p_1p_2p_3v}\, e^{\bar{i}(\tilde{\mathbf{p}}\bullet\bar{\boldsymbol{\omega}})t}\,|\mathbf{z}_r|^{\mathbf{m}}\, e_j \tag{54}$$

$$|\mathbf{z}_r|^{\mathbf{m}} = z_1^{m_1} z_2^{m_2}...z_r^{m_r},\ \bar{i}=\sqrt{-1},\ m_1+......+m_r = 2;\ \bar{\boldsymbol{\omega}} = \{\omega_1\ \ \omega_2\ \ \omega_f\};\ \tilde{\mathbf{p}}=\{p_1\ \ p_2\ \ p_3\}$$

As before, we can obtain $\bar{\mathbf{h}}_{23}(\mathbf{z}_r,t)$ via term-by-term comparison if and only if the following 'combined reducibility condition' is satisfied.

$$\bar{i}\,(\tilde{\mathbf{p}}\bullet\bar{\boldsymbol{\omega}}) + \sum_{l=1}^{r}(m_l\lambda_l) - \bar{\lambda}_p \neq 0 \tag{55}$$

It can be observed that this 'combined reducibility condition' involves a contribution from the forcing.

Once $\bar{\mathbf{h}}_{23}(\mathbf{z}_r,t)$ is known, we can solve the equation (47) which contains the contribution from $\mathbf{w}_{s3}(\mathbf{z}_r,t)$ (denoted by $\mathbf{w}_{s3_1}(\mathbf{z}_r,t)$) and $\mathbf{w}_{s2}(\mathbf{z}_r,t)$ (represented by $\mathbf{w}_{s3_1}(\mathbf{z}_r,t)$) and $\dfrac{\partial \bar{\mathbf{h}}_{23}(\mathbf{z}_r,t)}{\partial \mathbf{z}_r}\bar{\mathbf{F}}_r(t)$

To determine $\bar{\mathbf{h}}_{13}(\mathbf{z}_r,t)$, we expand the known and the unknown terms in the multiple Fourier series of the form given by equation (40) and equation (41), respectively and obtain the '*combined reducibility condition*' given by equation (43).

Further, to obtain $\overline{\mathbf{h}}_{03}(t)$, we use the convolution theorem, as before. However, it can be seen from equation (46) (which contains only temporal arguments), the forcing terms arise from $\dfrac{\partial \overline{\mathbf{h}}_{13}}{\partial \mathbf{z}_r}\overline{\mathbf{F}}_r(t)$ (forcing on the master states) and nonlinear terms $\mathbf{w}_{s2}(\mathbf{z}_r,t)$, $\mathbf{w}_{s3}(\mathbf{z}_r,t)$ and denoted as $\mathbf{w}_{s_q}(t)$, which can be expressed as

$$\mathbf{w}_{s_q}(t) = \sum_{j=1}^{s} \sum_{p_1=-\infty}^{+\infty} \sum_{p_2=-\infty}^{+\infty} \sum_{p_3=-\infty}^{+\infty} a_{j\overline{m}v} e^{\overline{i}(\tilde{\mathbf{p}}\bullet\overline{\boldsymbol{\omega}})t} \mathbf{e}_j \qquad (56)$$

where $\overline{i} = \sqrt{-1}$, $\overline{\boldsymbol{\omega}} = \{\omega_1 \ \omega_2 \ \omega_f\}$; $\tilde{\mathbf{p}} = \{p_1 \ p_2 \ p_3\}^T$.
The equation (46) cannot be solved if

$$\lambda_s = \tilde{\mathbf{p}} \bullet \overline{\boldsymbol{\omega}}. \qquad (57)$$

It is possible to continue the procedure discussed above to construct the relationship between 'slave' and 'master' states to the desired order and recover various '*resonance conditions*' involving contributions from external excitation and nonlinearity at multiple orders.

## 4. APPLICATIONS

Consider a coupled undamped Mathieu Hill-type nonlinear quasi-periodic system subjected to external excitation given by

$$\begin{aligned}\ddot{x} + (a_1 + b_1 \cos \omega_1 t + c_1 \cos \omega_2 t)x + x^2 y &= A_1 \cos(\omega t) \\ \ddot{y} + (a_2 + b_2 \cos \omega_1 t + c_2 \cos \omega_2 t))x + y^2 x &= A_2 \cos(\omega t)\end{aligned} \qquad (58)$$

The equation (58) can be expressed as

$$\dfrac{d}{dt}\begin{bmatrix}x \\ \dot{x} \\ y \\ \dot{y}\end{bmatrix} = \begin{bmatrix}\tilde{\mathbf{A}}_1(t) & \mathbf{0} \\ \mathbf{0} & \tilde{\mathbf{A}}_2(t)\end{bmatrix}\begin{bmatrix}x \\ \dot{x} \\ y \\ \dot{y}\end{bmatrix} + \begin{bmatrix}0 \\ x^2 y \\ 0 \\ y^2 x\end{bmatrix} + \begin{bmatrix}0 \\ \cos(\omega t) \\ 0 \\ \cos(\omega t)\end{bmatrix} \qquad (59)$$

where $\tilde{\mathbf{A}}_1(t) = \begin{bmatrix}0 & 1 \\ -(a_1 + b_1 \cos \omega_1 t + c_1 \cos \omega_2 t) & 0\end{bmatrix}$, $\tilde{\mathbf{A}}_2(t) = \begin{bmatrix}0 & 1 \\ -(a_2 + b_2 \cos \omega_1 t + c_2 \cos \omega_2 t) & 0\end{bmatrix}$

Applying the L-P transformation $\mathbf{x}(t) = \mathbf{Q}(t)\mathbf{z}(t)$ and its inverse to the equation (59) yields an equation similar to the equation (22)

$$\frac{d}{dt}\begin{bmatrix} z_1 \\ z_2 \\ z_3 \\ z_4 \end{bmatrix} = \begin{bmatrix} \mathbf{J}_1 & \mathbf{0} \\ \mathbf{0} & \mathbf{J}_2 \end{bmatrix}\begin{bmatrix} z_1 \\ z_2 \\ z_3 \\ z_4 \end{bmatrix} + \begin{bmatrix} f_{11}(\mathbf{z},t) \\ f_{12}(\mathbf{z},t) \\ f_{21}(\mathbf{z},t) \\ f_{22}(\mathbf{z},t) \end{bmatrix} + \begin{bmatrix} F_{11}(t) \\ F_{12}(t) \\ F_{21}(t) \\ F_{22}(t) \end{bmatrix} \qquad (60)$$

Where $\mathbf{J}_1 = \begin{bmatrix} \dfrac{\left(-\lambda_1 + \dfrac{b_1^2\left(C_2\lambda_3(-4\lambda_3+\lambda_3)+C_1\lambda_4(-4\lambda_1+\lambda_4)\right)}{(4\lambda_1-\lambda_3)\lambda_3(4\lambda_1-\lambda_4)\lambda_4}\right)}{\sqrt{\lambda_1}} & 0 \\ 0 & \dfrac{\left(\lambda_1 + b_1^2\left(\dfrac{C_1}{4\lambda_1\lambda_3-\lambda_3^2}+\dfrac{C_2}{4\lambda_1\lambda_4-\lambda_4^2}\right)\right)}{\sqrt{\lambda_1}} \end{bmatrix}$

$\mathbf{J}_2 = \begin{bmatrix} \dfrac{\left(-\lambda_2 + \dfrac{b_2^2\left(C_3\lambda_3(-4\lambda_3+\lambda_3)+C_4\lambda_4(-4\lambda_1+\lambda_4)\right)}{(4\lambda_2-\lambda_3)\lambda_3(4\lambda_2-\lambda_4)\lambda_4}\right)}{\sqrt{\lambda_2}} & 0 \\ 0 & \dfrac{\left(\lambda_2 + b_2^2\left(\dfrac{C_3}{4\lambda_2\lambda_3-\lambda_3^2}+\dfrac{C_4}{4\lambda_2\lambda_4-\lambda_4^2}\right)\right)}{\sqrt{\lambda_2}} \end{bmatrix}$

$\lambda_1 = a_1, \lambda_2 = a_2, \lambda_3 = \omega_1, \lambda_4 = \omega_2$ $C_1$ and $C_2$ are constant depending upon the initial conditions of the fictitious states. For more details on the computation of $\mathbf{J}$, we refer to reference [27]. At this point, we have to select master and slave states. Assuming eigenvalues of the $\mathbf{J}_1$ matrix are closer to $\omega$ (the frequency of external excitation), we choose $\mathbf{z}_r = \{z_1, z_2\}^T$ as the master states and $\mathbf{z}_s = \{z_3, z_4\}^T$ as the slave states and partition equation (60) similar to the equation (24).

$$\begin{Bmatrix} \dot{\mathbf{z}}_r \\ \dot{\mathbf{z}}_s \end{Bmatrix} = \begin{bmatrix} \mathbf{J}_r & 0 \\ 0 & \mathbf{J}_s \end{bmatrix}\begin{Bmatrix} \mathbf{z}_r \\ \mathbf{z}_s \end{Bmatrix} + \begin{Bmatrix} \mathbf{w}_r(\mathbf{z}_r,\mathbf{z}_s,t) \\ \mathbf{w}_s(\mathbf{z}_r,\mathbf{z}_s,t) \end{Bmatrix} + \begin{Bmatrix} \mathbf{F}_r(t) \\ \mathbf{F}_s(t) \end{Bmatrix} \quad \begin{matrix}(a)\\(b)\end{matrix} \qquad (61)$$

Where $\mathbf{J}_r = \mathbf{J}_1, \mathbf{J}_s = \mathbf{J}_2, \mathbf{w}_r(\mathbf{z}_r,\mathbf{z}_s,t) = \{f_{11}(\mathbf{z},t), f_{12}(\mathbf{z},t)\}^T, \mathbf{w}_s(\mathbf{z}_r,\mathbf{z}_s,t) = \{f_{21}(\mathbf{z},t), f_{22}(\mathbf{z},t)\}^T$
$\mathbf{F}_r(t) = \{F_{11}(t), F_{12}(t)\}^T, \mathbf{F}_s(t) = \{F_{21}(t), F_{22}(t)\}^T$

In this particular case,
$a_1 = 3, a_2 = 5, b_1 = b_2 = c_1 = c_2 = 2.5, \omega_1 = 2\pi\ rad/s, \omega_2 = 7\ rad/s, \omega = 1\ rad/s, A_1 = 1, A_2 = 1$
the $\mathbf{J}_1, \mathbf{J}_2, C_1$ and $C_2$ are given as

$$\mathbf{J}_1 = \begin{bmatrix} -1.78i & 0 \\ 0 & +1.78i \end{bmatrix}, \mathbf{J}_2 = \begin{bmatrix} -2.29i & 0 \\ 0 & +2.29i \end{bmatrix}, C_1 = \pi^2, C_2 = (3.5)^2 \qquad (62)$$

One can apply order reduction techniques discussed in section 3.

**a) Order reduction using the linear method**

Equation (61) comprises of 4 states $\{z_1 \quad z_2 \quad z_3 \quad z_4\}^T$. Following the procedure outlined in section 3.1, we neglect the contribution from the non-dominant states $\mathbf{z}_s = \{z_3 \quad z_4\}^T$, and the system dynamics is approximated by

$$\begin{Bmatrix} \dot{z}_1 \\ \dot{z}_2 \end{Bmatrix} = [\mathbf{J}_1] \begin{Bmatrix} z_1 \\ z_2 \end{Bmatrix} + \begin{Bmatrix} f_{11}(z_1, z_2, 0, 0, t) \\ f_{22}(z_1, z_2, 0, 0, t) \end{Bmatrix} + \begin{Bmatrix} F_{11}(t) \\ F_{12}(t) \end{Bmatrix} \tag{63}$$

The equation (63) is the reduced-order model of the system described by the equation (61). This reduced-order system is integrated numerically with typical initial conditions, and all the states in $\mathbf{x}$ are obtained using the L-P $\mathbf{x}(t) = \mathbf{Q}(t)\mathbf{z}(t)$ transformation. This solution is called linear reduced-order system response. This response can be compared with the original system's response calculated via numerical integration of the equation (58). The time trace comparison is shown in Figures 1 and 2. Figure 3 compares phase planes for the original and the reduced-order system via the linear technique. It can be noticed that the linear reduced-order model fails to capture the dynamics of the original system. One reason for this failure is that the slave states are also excited by forcing $\mathbf{F}_s(t)$ that is completely ignored in the reduced-order model. The linear order reduction technique may yield acceptable results when the eigenvalues corresponding to slave states have negative real parts or no forcing on slave states. However, in general, the linear order reduction approach for nonlinear quasi-periodic systems subjected to external excitation may not yield accurate results. For clarity, the Welch power spectrum for the original system response is compared with the Welch power spectrum for the linearly reduced system in Figure 4. These power spectrums do not match, indicating that the original system's dynamics (frequency content) are not captured in the linearly reduced-order system.

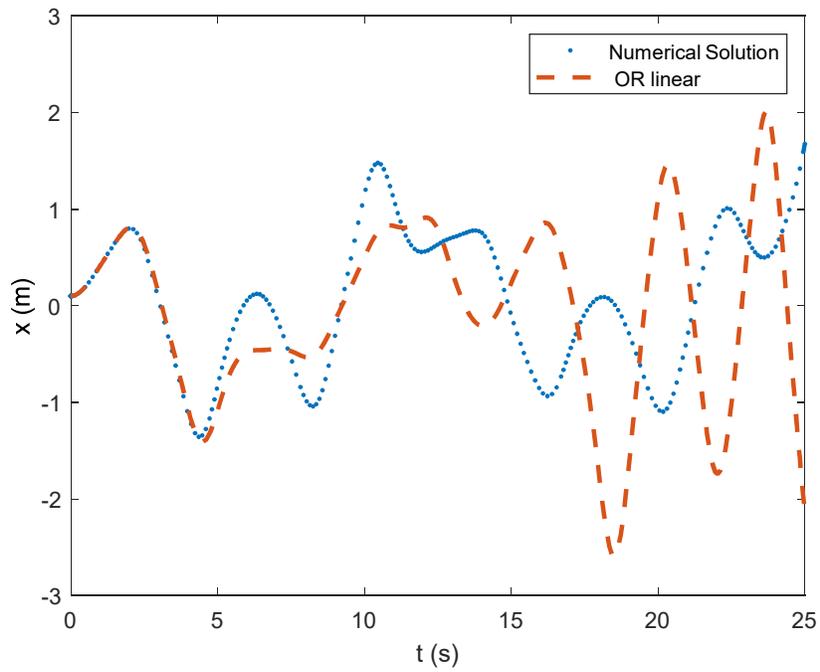

Figure 1: Time trace comparisons of original and the reduced-order system via the linear approach $x(t)$ v/s time

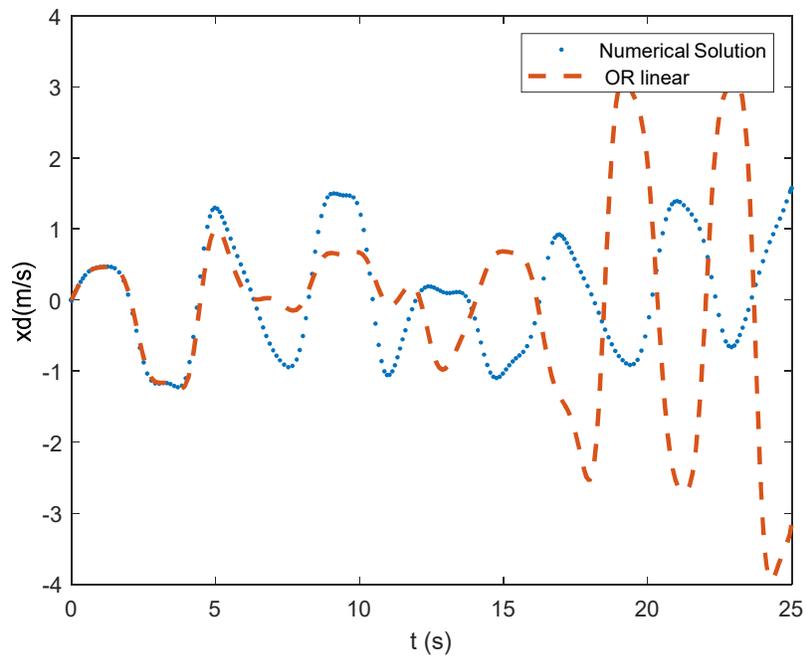

Figure 2: Time trace comparisons of original and the reduced-order system via the linear approach $\dot{x}(t)$ v/s time

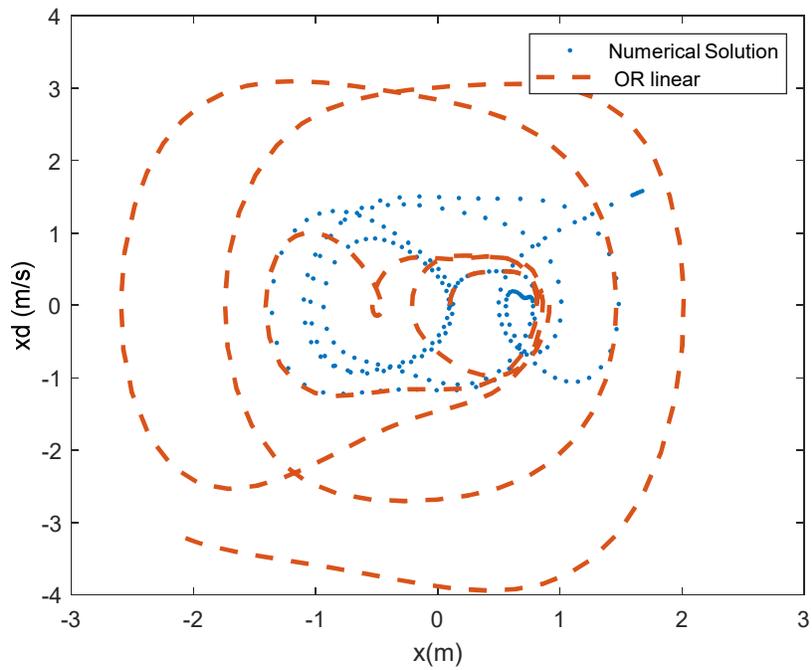

Figure 3: Phase plane comparisons of original and the reduced-order system via the linear approach $x(t)$ v/s $\dot{x}(t)$

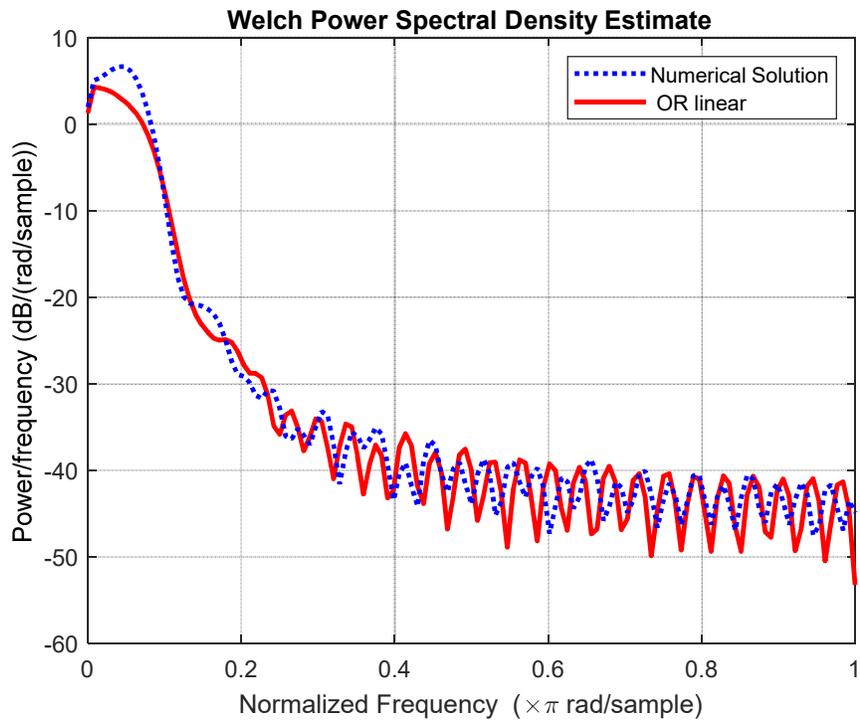

Figure 4: Welch Power Spectral density comparison of original and the reduced-order system via the linear approach

**b) Order reduction using an invariant manifold**

As discussed earlier, we try to relate the non-dominant states to the dominant states by a quasi-periodic nonlinear transformation. If the system does not exhibit any resonances (like the case under consideration), then the 'reducibility condition' is satisfied, and the system order can be reduced.

We start with the equation (61) and select the same states [$\mathbf{z}_r = \{z_1 \; z_2\}^T$] as the dominant states and try to find a nonlinear quasi-periodic relationship of the form given by the equation (27). For this particular example, the relationship between $\mathbf{z}_s$ and $\mathbf{z}_r$ are

$$\mathbf{z}_s = \sum_i \mathbf{h}_i(z_1, z_2, t) \equiv \mathbf{H}(z_1, z_2, t), \; s = 3, 4 \tag{64}$$

where $\quad \mathbf{h}_i = \sum_{\bar{\mathbf{m}}} \bar{\mathbf{h}}_i(t) z_1^{m_1} \ldots z_2^{m_2}, \quad \bar{\mathbf{m}} = (m_1, m_2)^T, \; m_1 + m_2 = 3 \tag{65}$

Here $\bar{\mathbf{h}}_i(t)$ are the unknown quasi-periodic vector coefficients. We substitute the equation (64) into the equation(61). After expanding $\bar{\mathbf{h}}_i(t)$ and $\mathbf{w}_s(\mathbf{z}_r, t)$ ($s = 3, 4$) in the Fourier series and neglecting the terms of higher-order, we obtain the relationship between the dominant and the non-dominant states as

$$z_3 = \mathbf{H}_1(z_1, z_2, t), \quad z_4 = \mathbf{H}_2(z_1, z_2, t) \tag{66}$$

The equation (66) is substituted into equation (61)-a to get the reduced-order model as

$$\begin{Bmatrix} \dot{z}_1 \\ \dot{z}_2 \end{Bmatrix} = [\mathbf{J}_1] \begin{Bmatrix} z_1 \\ z_2 \end{Bmatrix} + \begin{Bmatrix} \underline{w}_1(z_1, z_2, t) \\ \underline{w}_2(z_1, z_2, t) \end{Bmatrix} + \begin{Bmatrix} F_{11}(t) \\ F_{12}(t) \end{Bmatrix} \tag{67}$$

The equation (67) is the reduced-order model of the system described by the equation (61). As before, this reduced-order system is integrated numerically with typical initial conditions, and all the states in $\mathbf{x}$ are obtained using the L-P $\mathbf{x}(t) = \mathbf{Q}(t)\mathbf{z}(t)$ transformation. This solution is called nonlinear reduced-order system response. Similar to the linear reduced-order system analysis, the nonlinear reduced-order system response can be compared with the original system's response calculated via numerical integration of the equation (58). The time trace comparison is shown in Figures 5 and 6. Figure 7 compares phase planes for the original and the reduced-order system via the nonlinear technique. It can be noticed that the nonlinear reduced-order model captures the dynamics of the original system quite well. For additional insight, the Welch power spectrum for the original system response is compared with the Welch power spectrum for the nonlinearly reduced system in Figure 8. These power spectrums match, indicating that the original system's dynamics (frequency content) are captured in the nonlinearly reduced-order system. These symbolic computations were performed using Mathematica™.

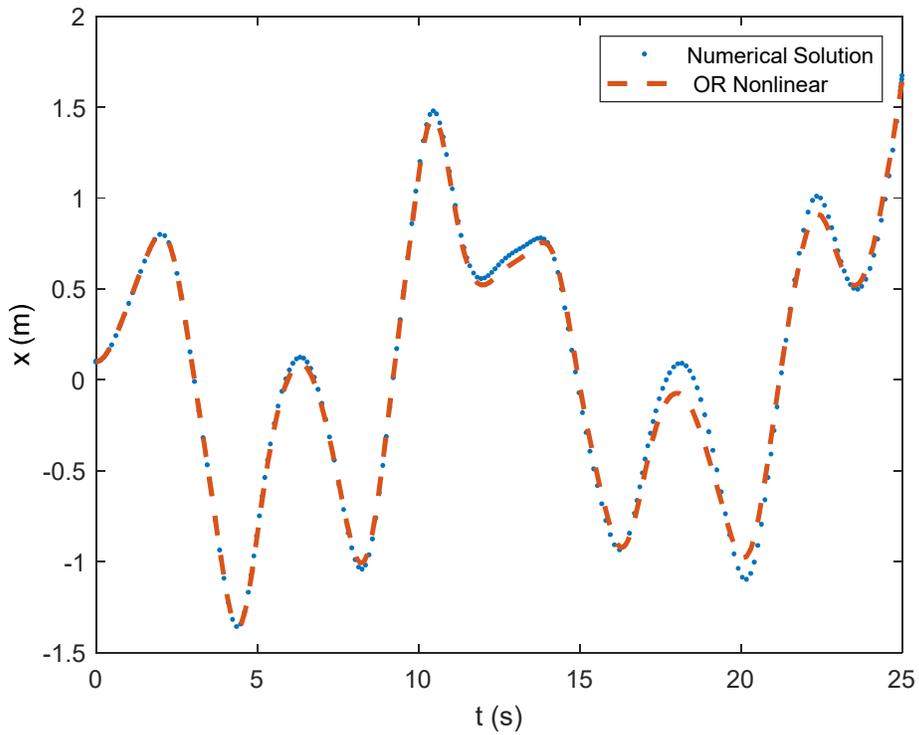

Figure 5: Time trace comparisons of original and the reduced-order system via the nonlinear approach $x(t)$ v/s time

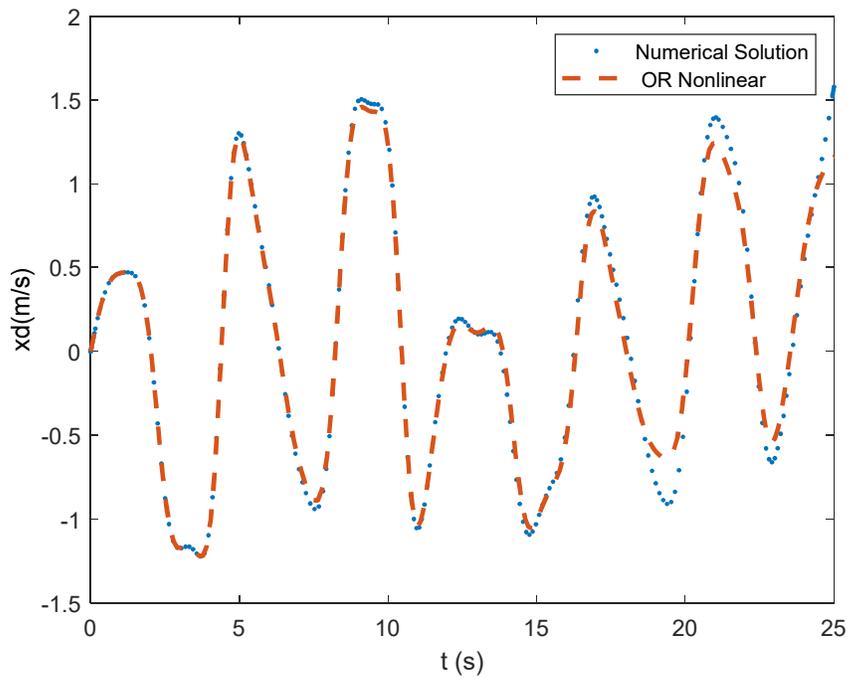

Figure 6: Time trace comparisons of original and the reduced-order system via the nonlinear approach $\dot{x}(t)$ v/s time

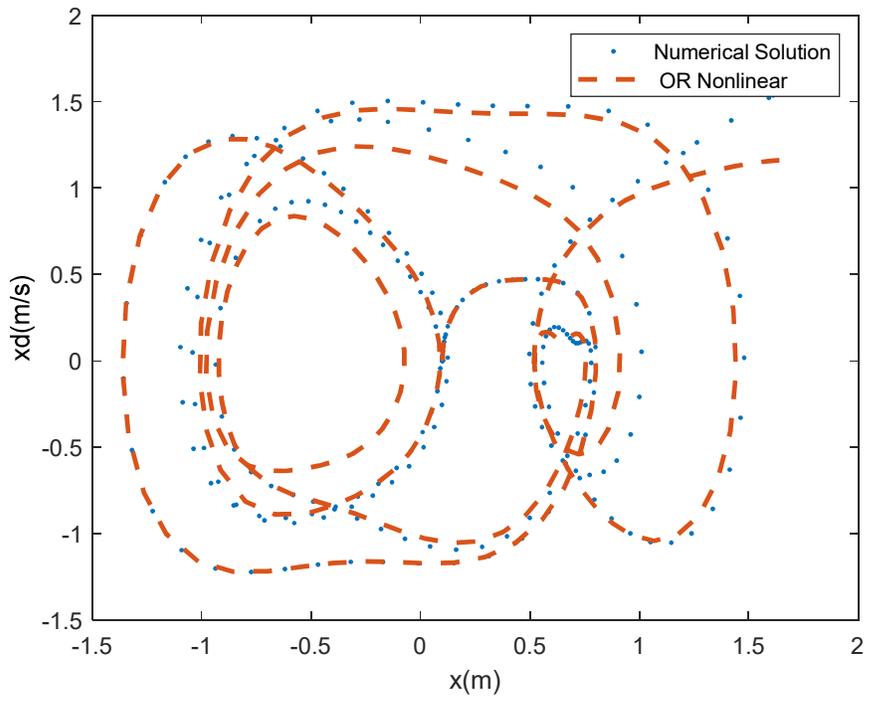

Figure 7: Phase plane comparisons of original and the reduced-order system via the nonlinear approach $x(t)$ v/s $\dot{x}(t)$

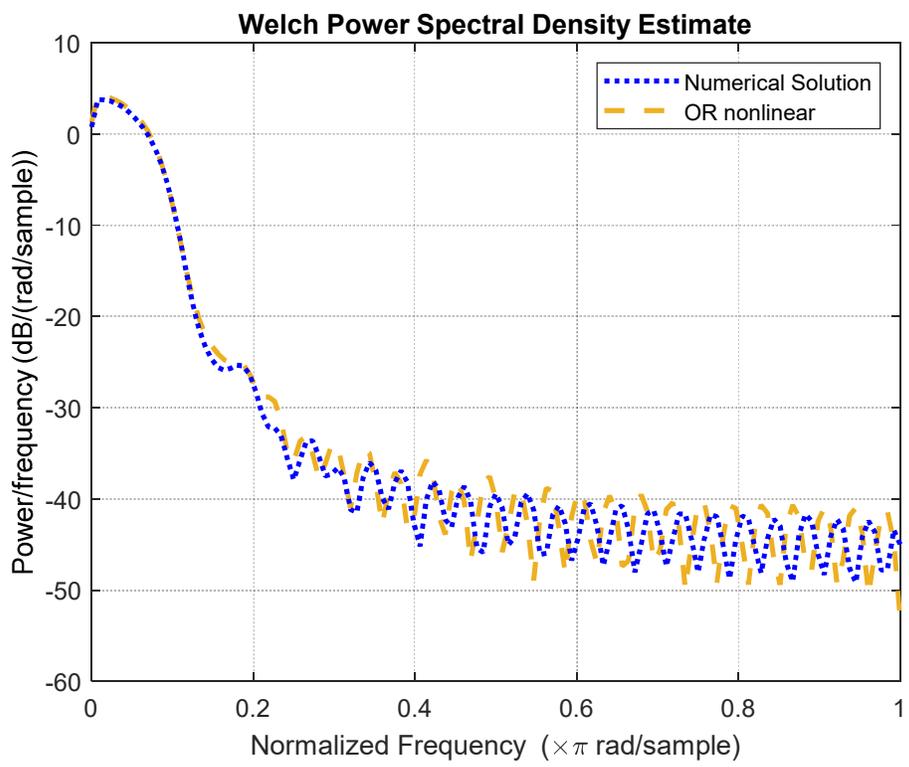

Figure 8: Welch Power Spectral density comparison of original and the reduced-order system via the nonlinear approach

## 5. DISCUSSION AND CONCLUSIONS

This paper presents a technique for obtaining a reduced-order model of a nonlinear quasi-periodic system subjected to external excitation. The central idea here is to assume a quasi-periodic transformation with unknown coefficients between the master and the slave states. This transformation can be determined by collecting the terms of the same order and solving them using harmonic balance. In the solution process, we obtain reducibility conditions that indicate resonances between system states, nonlinearity, and external excitation. The linear resonance condition is also obtained as we find the solution of quasi-periodic transformation. One crucial point here is how one can decide which states to retain and which ones to eliminate. Initially, one could start with the states corresponding to eigenvalues close to external excitation frequencies and start the order reduction process. The resonance interactions in the nonlinear quasi-periodic system are complex, and in the course of order reduction, one may see a "small deviser" problem. Such a case indicates resonant interaction, and these resonant states must be included in the master states. It can be noted that with the advent of symbolic software like Mathematica and Maple, the procedure for order reduction can be automated [41]. One can consider quasi-periodic and external excitation as fictitious states and carry out the order reduction. This approach is presented in reference [42], and further simplification via the method of form can be achieved as discussed in reference [43] for autonomous systems.

The reduced-order model will contain all the essential dynamics and responses of the reduced-order system quantitatively and qualitatively, similar to the original system. This reduced-order system can be simplified using the method of normal forms. One can study this simplified system for bifurcation and control. The reduced-order system can be used for the optimization of essential parameters, study sensitivity, and design controllers.

**FUNDING:** This work was not supported by any funding agency.

## NOMENCLATURE

$\mathbf{x}$ - $n$ vector of states
$\mathbf{A}(t)$ - $n \times n$ time quasi-periodic matrix
$\mathbf{f}(\mathbf{x},t)$ - nonlinear $n$ vector such that $\mathbf{f}(0,t) = 0$
$\mathbf{Q}(t)$ - L-P transformation matrix of dimension $n \times n$
$\mathbf{M}$ - Modal matrix of dimension $n \times n$
$\mathbf{z}$ - $n$ vector of the L-P transformed states
$\mathbf{z}_r$ - $r(r \ll n)$ vector of dominant states
$\mathbf{z}_s(t)$ - $s(s+r=n)$ vector of non-dominant (slave) states
$\mathbf{J}_r$ - $r \times r$ Jordan block corresponding to dominant states
$\mathbf{H}(\mathbf{z}_r,t)$ - Nonlinear quasi-periodic invariant manifold function relating the non-dominant states to dominant states